\documentclass[12pt]{article}
\usepackage{cite}
\usepackage[margin=2.5cm]{geometry}
\usepackage{amsmath,amssymb,amsfonts}
\usepackage{graphicx}
\usepackage{textcomp}
\usepackage{amsmath,amssymb,bm,bbm,mathrsfs,amscd}
\usepackage{calc}
\usepackage{color}
\usepackage{amsthm}
\usepackage{dsfont}
\usepackage{graphicx}
\usepackage{pifont}
\usepackage{amsmath}
\usepackage{amsfonts}
\usepackage{amssymb}
\usepackage{rotating}
\usepackage{mathtools}
\renewcommand\footnotemark{}
\usepackage{color}
\usepackage[caption=false]{subfig}
\usepackage{enumerate}
\newtheorem{theorem}{Theorem}
\newtheorem{problem}{Problem}
\def\diag{{\rm diag}\,}
\renewcommand{\eqref}[1]{Eq.~(\ref{#1})}  
\newcommand{\vect}[1]{\boldsymbol{#1}}
\def\BibTeX{{\rm B\kern-.05em{\sc i\kern-.025em b}\kern-.08em
    T\kern-.1667em\lower.7ex\hbox{E}\kern-.125emX}}

\begin{document}
\title{Analysis, Prediction, and Control of Epidemics: A Survey from Scalar to Dynamic Network Models}
\author{Lorenzo Zino and Ming Cao\thanks{© 2021 IEEE.  Personal use of this material is permitted.  Permission from IEEE must be obtained for all other uses, in any current or future media, including reprinting/republishing this material for advertising or promotional purposes, creating new collective works, for resale or redistribution to servers or lists, or reuse of any copyrighted component of this work in other works.}}
\date{\normalsize Faculty of Science and Engineering, University of Groningen, 9747 AG Groningen, Netherlands (\texttt{\{lorenzo.zino, m.cao\}@rug.nl})}

\maketitle

\begin{abstract}
During the ongoing COVID-19 pandemic, mathematical models of epidemic spreading have emerged as powerful tools to produce valuable predictions of the evolution of the pandemic, helping public health authorities decide which intervention policies should be implemented. The study of these models --- grounded in the systems theory and often analyzed using control-theoretic tools --- is an extremely important research area for many researchers from different fields, including epidemiology, engineering, physics, mathematics, computer science, sociology, economics, and management. In this survey, we review the history and present the state of the art in the modeling, analysis, and control of epidemic dynamics. We discuss different approaches to epidemic modeling, either deterministic or stochastic, ranging from the first implementations of scalar systems of differential equations to describing the epidemic spreading at the population level, and to more recent models on dynamic networks, which capture the spatial spread and the time-varying nature of human interactions.
\end{abstract}

\section{Brief History of 260 years of Mathematical Models of Epidemics}

Since the beginning of human history, pandemics have posed deadly threats, which often decimated our species. Hence, it was not surprising that, in parallel with the theoretical development of calculus, mathematicians started to apply their theoretical paradigms to describe, study, and unveil the mechanisms of spreading of infectious diseases. In this vein, the first milestone can be found in the work on smallpox by Daniel Bernoulli~\cite{Bernoulli1760}, published in 1760. In the 18th century, there was an ongoing public debate about variolation, i.e., inoculation of infectious material from smallpox cases to induce a mild infection and lifelong immunity. In Bernoulli's work, the Dutch-born Swiss mathematician built a mathematical theory to support the effectiveness of variolation, helping its increasing adoption until the development of the smallpox vaccine in 1796. A latter milestone of mathematical modeling of epidemics is materialized in the studies of the 1849 and 1854 Cholera outbreaks in London by William Farr~\cite{Farr} and John Snow~\cite{Snow}, respectively.

However, it is not until the beginning of the 20th century that  differential and difference equations started being adopted as tools to model and analyze the spread of epidemic diseases. The very origin of this approach can be found in a paper by William Heaton Hamer~\cite{Hamer1906}, in which the British epidemiologist pioneered the use of a nonlinear formula to model the rate of the contagion process, proportional to the product between the number of susceptible individuals and the number of infectious individuals in the population. Such a modeling approach has been formalized and popularized by William Ogilvy Kermack and Anderson Gray McKendrick. In their seminal paper~\cite{Kermack1927}, and in two subsequent works~\cite{Kermack1932,Kermack1933}, the two Scottish researchers laid out the basis for the current mathematical theory of epidemic modeling, namely the susceptible--infected--susceptible (SIS) and the susceptible--infected--removed (SIR) models and the concept of epidemic threshold (and consequently, phase transition).

\subsection{Population SIS model}

In the SIS model, it is assumed that infected individuals after recovering do not acquire immunity, and thus become again susceptible to the disease. This model became very popular to study recurrent and endemic diseases, including sexually transmitted diseases, such as Gonorrhea~\cite{Yorke1978}. Here, we present the deterministic SIS model in a continuous-time framework, i.e., employing differential equations (as in the original contribution by Kermack and McKendrick~\cite{Kermack1927}). Note that, discrete-time implementations of the model using difference equations have also been extensively studied. More details on discrete-time models can be found, in e.g.,~\cite{Allen1994}. 

\begin{figure}
    \centering
        \subfloat[]{\includegraphics[scale=1.8]{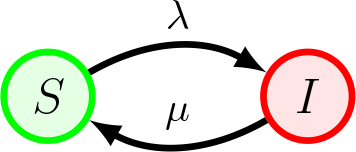} \label{fig:SIS}}\qquad
    \subfloat[]{\includegraphics[scale=1.8]{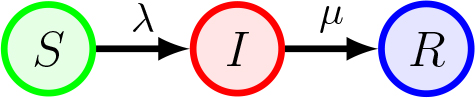}  \label{fig:SIR}}
    \caption{Schematics of the (a) SIS and (b) SIR models.}
    \label{fig:models}
\end{figure}

Formally, a unit mass population is split between two compartments: susceptible and infected. Let $s(t)\in[0,1]$ and $x(t)\in[0,1]$ be the fraction of susceptible and infected individuals in the population at time $t\geq 0$, respectively. The health state of the population evolves according to two mechanisms: contagion and recovery, illustrated in Fig.~\ref{fig:SIS}. Susceptible individuals may become infected if they interact with infected individuals, while infected individuals spontaneously recover, becoming again susceptible to the disease. According to the contagion mechanism, the total number of newly infected individuals is a nonlinear expression proportional to the number of susceptible individuals and to the number of infectious individuals in the population (as proposed in\cite{Hamer1906}), i.e., equal to \begin{equation}\lambda s(t)x(t)\,,\end{equation} where $\lambda>0$ is a positive parameter termed \emph{contagion rate}. The total number of recoveries in the population is proportional to the number of infected individuals, yielding the term $\mu x(t)$, where $\mu>0$ is the \emph{recovery rate}.

\begin{figure}
\centering
\subfloat[]{\includegraphics[scale=1.8]{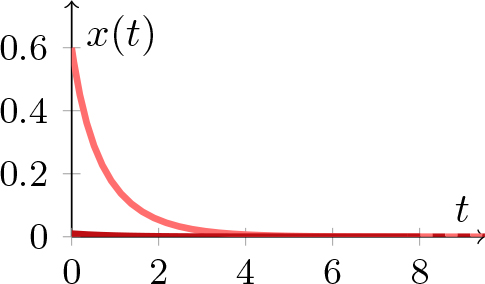}}\qquad
\subfloat[]{\includegraphics[scale=1.8]{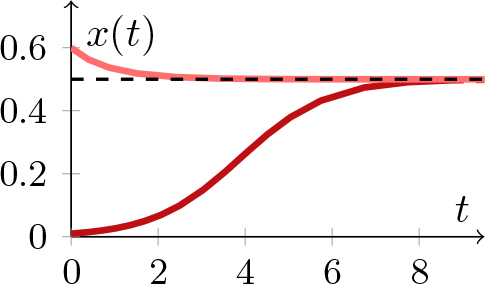}}
    \caption{Two trajectories with different initial conditions of the population SIS model in~\eqref{eq:sis} with (a) $\lambda/\mu=0.5$ (i.e., below the epidemic threshold) and (b) $\lambda/\mu=2$ (i.e., above the epidemic threshold). Below the epidemic threshold, both trajectories converge to the disease free equilibrium; above the epidemic threshold, they both converge to the endemic equilibrium $x^*=0.5$, denoted by the black dashed line.}
    \label{fig:SIS_simulations}
\end{figure}

Hence, the heath state of the population is described by the two-dimensional state variable $[s(t),\,x(t)]^\top$, which evolves according to the following system of ordinary differential equations (ODEs):
\begin{equation}\label{eq:sis_system}
\left\{\begin{array}{lll}
    \dot s(t)&=&-\lambda s(t)x(t)+\mu x(t)\\[6pt]
    \dot x(t)&=&\lambda s(t)x(t)-\mu x(t)\,.
\end{array}\right.
\end{equation}
Note that the two equations are linearly dependent, since the total mass of the population is preserved, that is, $\dot s(t)+\dot x(t)=0$. Hence, the system of ODEs in \eqref{eq:sis_system} can be reduced to the single nonlinear ODE:
\begin{equation}\label{eq:sis}
    \dot x(t)=\lambda (1-x(t))x(t)-\mu x(t)\,,
\end{equation}
from which one can immediately observe that the domain $[0,1]$ is positively invariant. From the analysis of this equation and its equilibrium points, one can conclude that, depending on the model parameters $\lambda$ and $\mu$, two outcomes may occur: either the origin (\emph{disease-free equilibrium}) is the unique equilibrium of the system in $[0,1]$ and it is globally asymptotically stable; or the origin becomes unstable, giving rise to a (almost) globally asymptotically stable \emph{endemic equilibrium} $x^*\in(0,1]$. These two behaviors are illustrated in Fig.~\ref{fig:SIS_simulations}. The phase transition between these two behaviors occur when a certain relation between the model parameters is satisfied, which is called \emph{epidemic threshold}. In this survey, we opt for expressing such a threshold in terms of the ratio between the contagion rate $\lambda$ and the recovery rate $\mu$. The following result from~\cite{Kermack1927} fully characterizes the behavior of the SIS model.

\begin{theorem}\label{th:sis}
If the population SIS model in \eqref{eq:sis} satisfies ${\lambda}/{\mu}\leq1$, then the disease free equilibrium $x=0$ is globally asymptotically stable. Otherwise, if ${\lambda}/{\mu}>1$, the disease-free equilibrium is unstable and \eqref{eq:sis} has a (almost) globally asymptotically stable endemic equilibrium, corresponding to $x^*=1-\mu/\lambda$.
\end{theorem}

The epidemic threshold is often associated with the \emph{basic reproduction number} $\mathcal R_0$, that is, the average number of contagions that a single infected person will cause. The concept of basic reproduction number, although already touched upon in the original works on the SIS model~\cite{Kermack1927}, was not formally introduced until the work by George MacDonald in the early 1950s~\cite{macdonald1952}. For the SIS model in \eqref{eq:sis}, in fact, $\mathcal R_0=\lambda/\mu$. However, in more complex models, e.g., those including additional contagion mechanisms, stochasticity, heterogeneity, and  dynamic network structures, the basic reproduction number may not be sufficient to characterize the behavior of the system, which may depend on other factors such as its variability~\cite{Li2011}. For this reason, in this survey we prefer to present the results in terms of the epidemic thresholds.

\subsection{Extensions and limitations of population epidemic models}

\begin{figure}
\centering
\subfloat[]{\includegraphics[scale=1.8]{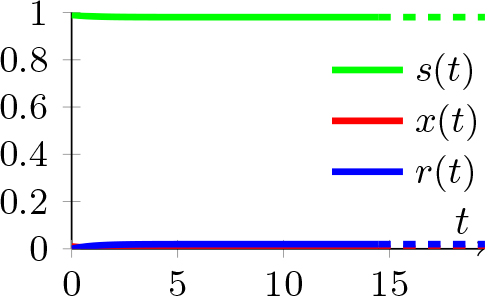}}\qquad
\subfloat[]{\includegraphics[scale=1.8]{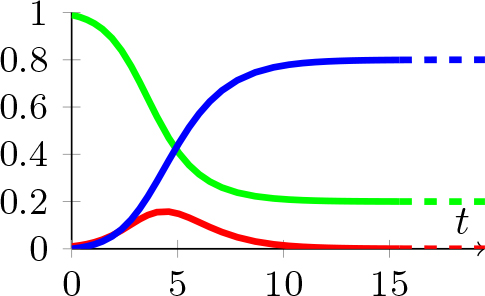}}
    \caption{Trajectories of the population SIR model in~\eqref{eq:sir} with (a) $\lambda/\mu=0.5$ (i.e., below the epidemic threshold) and (b) $\lambda/\mu=2$ (i.e., above the epidemic threshold). The green, red, and blue curves represent the fraction of susceptible $s(t)$, infected $x(t)$, and removed $r(t)$ individuals, respectively. Below the epidemic threshold, the disease is quickly eradicated and almost all of the population remains susceptible; above the epidemic threshold, one observes an initial growth of the epidemic prevalence and the majority of the population is eventually removed.
    }
    \label{fig:SIR_simulations}
\end{figure}

In the population SIR model (see Fig.~\ref{fig:SIR}), which was proposed by W.~O. Kermack and A.~G. McKendrick in their seminal works with the SIS model, a further compartment named \emph{removed} is introduced to account for individuals that recover from the disease and become immune, or die; $r(t)\in[0,1]$ is the fraction of population in the removed state. Hence, the heath state of the system is defined by the three-dimensional state $[s(t),\,x(t)\,,r(t)]^\top$ (with $s(t)+x(t)+r(t)=1$), which evolves according to the following system of ODEs:
\begin{equation}\label{eq:sir}
\left\{\begin{array}{lll}
    \dot s(t)&=&-\lambda s(t)x(t)\\[6pt]
    \dot x(t)&=&\lambda s(t)x(t)-\mu x(t)\\[6pt]
    \dot r(t)&=&\mu x(t)\,,
\end{array}\right.
\end{equation}
where the three equations are linearly dependent, since $\dot s(t)+\dot x(t)+\dot r(t)=0$. Hence, the system of ODEs in \eqref{eq:sir} reduces to a planar system. The analysis of such a planar system (see, e.g.,~\cite{Mei2017}) shows that the ratio $\lambda/\mu$ also plays an important role in the SIR model, modulated by the initial fraction of susceptible individuals $s(0)$. If $\lambda s(0)/\mu<1$, then the fraction of infected individuals exponentially converges to $0$. If $\lambda s(0)/\mu>1$, then the outbreak monotonically grows until reaching an epidemic peak and, only after the peak is reached, decreasing to $0$. Two trajectories of the population SIR model to illustrate these two different behaviors are shown in Fig.~\ref{fig:SIR_simulations}.

Toward the development of more realistic models, further compartments have been introduced. This allows to account for many realistic phenomena, including latency periods between contagion and infectiveness, multiple stages of the progression of the disease, as well as to model pharmaceutical and nonpharmaceutical interventions, including vaccination and quarantine. More details on this extensive family of mathematical models can be found in~\cite{Hethcote2000,brauer2011mathematical}. This large variety of models has been widely adopted in the recent years for many theoretical and practical studies of epidemic processes, including seasonal and pandemic influenza~\cite{Coburn2009}, SARS~\cite{Ng2003}, Ebola~\cite{Legrand2007}, and the ongoing COVID-19 outbreak~\cite{Giordano2020, Zhong2020, Casella2020}.

Control of epidemics was one of the major motivations for the development of the first mathematical models of epidemics. In fact, Bernoulli proposed his model to support the use of variolation~\cite{Bernoulli1760}. The idea of exploiting the dynamical system formulation of epidemic models to derive and assess control techniques can be traced back to the work on Malaria by Robert Ross~\cite{Ross1910}, in which the British medical doctor claimed that the reduction of the population of mosquitoes below a certain threshold would suffice to eradicate the disease. The development of compartmental models has laid the foundation of a rigorous mathematical framework, grounded in the theory of dynamical systems, to study the spread of epidemic processes and assess the effectiveness of different control strategies implemented to mitigate and combat them, by employing control-theoretic techniques.

However, this class of deterministic compartmental models suffers from some inherent limitations, which may hinder their applicability to real-world scenarios, thus calling for several extensions to be discussed in the rest of this survey. The first limitation is concerned with the use of differential equations to model the dynamical process. This usage is under the assumption of having a large population that can be approximated by a continuous distribution; it also assumes that the presence of noise and uncertainties can be omitted. However, epidemics often start from small local outbreaks, where finite population effects and stochasticity cannot be neglected. To account for these important factors, compartmental models of epidemics in populations have been extended to stochastic frameworks, by using the theory of Markov processes and branching processes. Exhaustive surveys of these models can be found in the following books~\cite{Bailey1975,Hethcote2000}.

The second limitation lies in the fact that this class of deterministic compartmental models provides a description of the epidemic process only at the population level, that is, in terms of the overall fractions of susceptible and infected individuals in the system. In real-world epidemic outbreaks, the pattern of interactions between the members of a population --- and, consequently, the contagion process --- has a complex structure, which may be influenced by the spatial location of the individuals of a population and the pattern of their interactions, calling for a finer representation of the population structure. In the last few decades, network theory has emerged as a powerful tool to capture such a complexity and formalize it in mathematically tractable models~\cite{keeling2011modeling,RevModPhys.87.925,Nowzari2016,Mei2017}. 

In the rest of this survey, we will focus on network epidemic models, mostly discussing the SIS model, which is among the simplest and most studied models. In Section \ref{sec:static}, we will review some classical results for its deterministic and stochastic implementations on static networks. In Section \ref{sec:dynamic}, we will move to the most recent developments and the state of the art for the research on epidemics on dynamic networks. Section \ref{sec:control} is devoted to the discussion of control of deterministic and stochastic network epidemic models, with specific attention on the control of dynamic networks. In Section \ref{sec:covid}, we present a focused discussion on the current challenge of the ongoing COVID-19 pandemic and on how the mathematical modeling of epidemics can provide powerful tools toward fighting against and mitigating its spreading. Finally, in Section \ref{sec:future}, we briefly summarize the content of this survey and outline some promising challenges for the future research.

\section{Classic Models of Epidemics on Networks}\label{sec:static}

Here, we present the network SIS model and report the main results for its deterministic or stochastic implementations. Then, we briefly mention some important extensions of these results to scenarios that incorporate further features of real-world systems. Towards this end, we first gather some definitions and basic notions on network and graph theories, which will be used throughout the survey.

\subsection{Basic notions on networks and graph theory}

\begin{figure}
\centering
\includegraphics[scale=1.8]{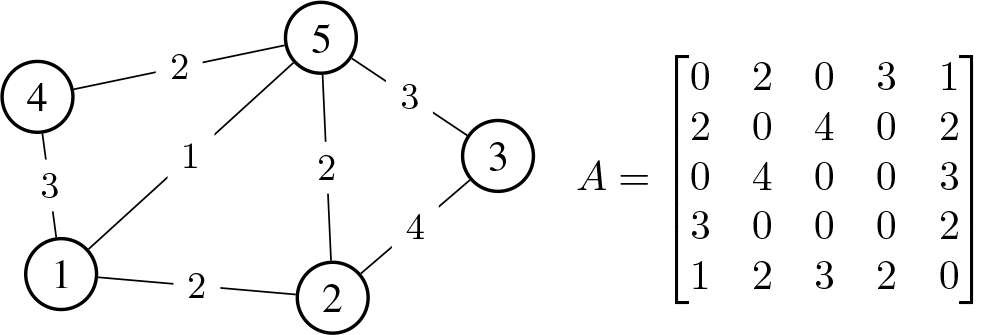}
\caption{Example of a weighted undirected static graph with $n=5$ nodes and its weighted adjacency matrix $A$. Note that the graph is connected.}
\label{fig:graph}
\end{figure}

A static network is represented by means of a time-invariant graph with $n$ nodes, denoted by a set of positive integer indices $\mathcal V=\{1,\dots,n\}$. Nodes are connected through a set of directed links $\mathcal E\subseteq \mathcal V\times\mathcal V$, such that $(i,j)\in\mathcal E$ if and only if node $i$ is connected to node $j$. The intensity of such connections is measured by a nonnegative quantity. Specifically, for any pair of nodes $i,j\in\mathcal V$, we define the \emph{connection strength} $a_{ij}\geq 0$ that quantifies how strongly node $i$ is connected to node $j$, and so, $a_{ij}>0\iff(i,j)\in\mathcal E$. The connection strengths are collected in a \emph{(weighted) adjacency matrix} $A\in\mathbb{R}_{\geq 0}^{n\times n}$. The triple $\mathcal G=(\mathcal V,\mathcal E, A)$ defines the \emph{graph}. An example is illustrated in Fig.~\ref{fig:graph}. 

 A network is \emph{undirected} if the corresponding adjacency matrix is symmetric, that is, $A=A^\top$; otherwise it is said to be \emph{directed}. A network is \emph{connected} (\emph{strongly connected}, for directed networks) if its adjacency matrix $A$ is irreducible, that is, if for any pair of nodes $i$ and $j$, there exists a sequence of nodes $v_1=i,v_2,\dots,v_k=j$ such that $(v_\ell,v_{\ell+1})\in\mathcal E$, for $\ell=1,\dots,k-1$. A network is \emph{unweighted} if the adjacency matrix of the corresponding graph $A$ has binary entries, that is, all nonzero entries (corresponding to links) are equal to $1$. Given a node $i\in\mathcal V$, we denote by $k_i=\sum_{j}a_{ij}$ its (weighted) \emph{degree}. Note that if the network is unweighted, then the degree $k_i$ is equal to the number of nodes that node $i$ is connected to, which are called \emph{neighbors} of $i$.

A dynamic network is represented by means of a time-varying graph $\mathcal G(t)=(\mathcal V,\mathcal E(t), A(t))$, where the time $t$ can be a discrete or a continuous index. The $n$ nodes in the node set $\mathcal V=\{1,\dots,n\}$ are time-invariant and connected through a time-varying set of links $\mathcal E(t)$. The matrix $A(t)\in\mathbb R_{\geq 0}^{n\times n}$ is the time-varying (weighted) adjacency matrix and measures the strengths of the connections between nodes at time $t$.

\subsection{Deterministic network models}

From the first implementation of a deterministic SIS model on a (static) network, proposed by Ana Lajmanovich and James A. Yorke in 1976 to study the spread of gonorrhea in a population~\cite{Lajmanovich1976}, network epidemic models have become an established and successful paradigm to study the spread of epidemic diseases in complex populations~\cite{keeling2011modeling, RevModPhys.87.925, Nowzari2016, Mei2017}. 

\begin{figure}
\centering
\includegraphics[scale=1.8]{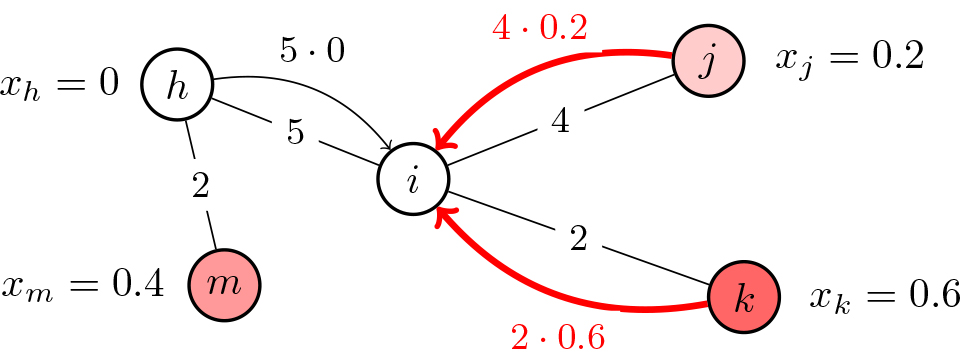}
\caption{Example of the contagion rate in a node of the network. The contagion rate in node $i$ is equal to the sum of the contribution coming from all its neighbors. Specifically, nodes $j$ and $k$ contribute to the rate, proportionally to the fraction of infected individuals in the nodes and the connection strengths; nodes $h$ and $m$ do not contribute: in $h$ there are no infected individuals, while $m$ is not a neighbor of $i$.}
\label{fig:network_sis}
\end{figure}

In the network SIS model, each node represents a group of individuals and is characterized by its health state, that is, ($s_i(t)$, $x_i(t)$). The health state of node $i\in\mathcal V$ represents the fraction of individuals per health states in the $i$th group~\cite{Lajmanovich1976}. Hence, $s_i(t)\in[0,1]$, $x_i(t)\in[0,1]$, and, similar to the population SIS model described in the previous section, it holds true that $x_i(t)+s_i(t)=1$. Hence, the health state of each node can be fully determined by the unique variable $x_i(t)$ and the entire state of the population is thus represented by the $n$-dimensional vector $\vect x(t)\in[0,1]^n$. Similar to the scalar model, the contagion rate in node $i$ has a nonlinear expression, equal to the product of the infection rate $\lambda_i$, the fraction of susceptible individuals in the $i$th group $s_i(t)$ and the strength of the interactions with other infected individuals in the population $m_i(t)$, which depends on the network structure and is equal to\begin{equation}\label{eq:contagion}m_i(t)=\sum_{j\in\mathcal V}a_{ij}x_j(t)\,,\end{equation} 
as illustrated in the simple example in Fig.~\ref{fig:network_sis}. For the sake of simplicity, in this section, we will present our main results under the assumption that the graph $\mathcal G=(\mathcal V,\mathcal E, A)$ is strongly connected. Results for directed networks that are not strongly connected can be found in~\cite{Meyers2006,Khanafer2016}.

Hence, the network SIS model is governed by a set of $n$ ODEs, one for each node of the network, where the equation that determines the evolution of $x_i(t)$, for a generic node $i\in\mathcal V$, is the following: 
\begin{equation}\label{eq:sis_network}
    \dot x_i(t)=\lambda_i (1-x_i(t))\sum_{j\in\mathcal V}a_{ij}x_j(t)-\mu_i x_i(t)\,,
\end{equation}
where $\lambda_i$ and $\mu_i$ are the infection and recovery rate of node $i$, respectively. While in most of the literature it is assumed that the epidemic parameters are homogeneous in the population, that is, $\lambda_i=\lambda$ and $\mu_i=\mu$, for all $i\in\mathcal V$, the general formulation in \eqref{eq:sis_network} with heterogeneous parameter has been adopted to model and study different scenarios, for instance, to model diseases that affect different age groups differently, which is key to implement and study targeted interventions and vaccination policies.

Note that, after introducing the $n$-dimensional vectors $\vect\lambda$ and $\vect\mu$ to gather all the $\lambda_i$ and $\mu_i$, $i\in\mathcal V$, respectively, the dynamics can be written in a compact, vector form as
\begin{equation}\label{eq:sis_network_matrix}
    \dot{\vect x}(t)=\Lambda\diag(\vect 1-\vect x(t))A\vect x(t)-M \vect x(t)\,,
\end{equation}
where $\Lambda=\diag(\vect\lambda)$, $M=\diag(\vect \mu)$, and $\vect 1$ is an $n$-dimensional all-$1$ vector from which is straightforward to observe that the domain $[0,1]^n$ is positively invariant. The following result, initially presented in~\cite{Lajmanovich1976}, extends the results for the population SIS model in Theorem \ref{th:sis} to the network SIS model. Different techniques have been used to prove this extension, including Lyapunov arguments~\cite{Fall2009} and positive systems theory~\cite{Khanafer2016}.
\begin{theorem}\label{th:sis_network}
Consider the homogeneous network SIS model in \eqref{eq:sis_network} with $\lambda_i=\lambda$ and $\mu_i=\mu$, for all $i\in\mathcal V$, on a strongly connected network. If  \begin{equation}\label{eq:r0_det}
     \frac{\lambda}{\mu}\leq\frac{1}{\rho(A)}\,,
\end{equation} where $\rho(A)$ is the spectral radius of the (weighted) adjacency matrix $A$, then the disease free equilibrium $\vect x=\vect 0$ is globally asymptotically stable. Otherwise, if $\lambda/\mu>1/\rho(A)$, the disease-free equilibrium is unstable and \eqref{eq:sis_network_matrix} has a unique (almost) globally asymptotically stable endemic equilibrium $\vect x^*$.
\end{theorem}
Different from the population SIS model, where the endemic equilibrium $\vect x^*$ has a closed-form expression depending on the model parameters, for the network SIS model such an explicit formula cannot be derived in general. However, iterative algorithms to compute such an equilibrium have been proposed in the literature. See, e.g.,~\cite{vanmieghem2009,Mei2017}. 

The convergence result in Theorem \ref{th:sis_network} can be easily extended to heterogeneous SIS models, where nodes have different contagion rates $\lambda_i$ and/or recovery rates $\mu_i$, as shown in~\cite{Lajmanovich1976, Wang2003, Khanafer2016}. In the heterogeneous case, the behavior of the system is determined by the spectral radius of the matrix $\Lambda A-M$. Specifically, if
\begin{equation}\label{eq:threshold_het}
   \rho(\Lambda A-M)\leq 1\,,
\end{equation}
the disease free equilibrium $\vect x=\vect 0$ is globally asymptotically stable. Otherwise, the SIS dynamics in \eqref{eq:sis_network_matrix} converge to the unique (almost) globally asymptotically stable endemic equilibrium $\vect x^*$. Note that this expression reduces to \eqref{eq:r0_det} in the homogeneous scenario. A technique to approximate the solution of \eqref{eq:sis_network_matrix} about the epidemic thresholds has been proposed in~\cite{Prasse2020}. 

Similarly, the SIR model and more complex compartmental models have been embedded and studied on network structures. For more details on these implementation and their analysis, we refer to the following review papers~\cite{Nowzari2016,Mei2017,RevModPhys.87.925}.

\subsection{Stochastic network models}

In their stochastic implementation, network epidemic models are defined as follows. Each node of the network represents an individual and is characterized by the health state $X_i(t)$, which coincides with one of the compartments. In the stochastic network SIS model, we have
\begin{equation}
X_i(t)=\left\{\begin{array}{ll}0&\text{if }i\text{ is susceptible at time }t\,,\\1&\text{if }i\text{ is infected at time }t\,.\end{array}\right.
\end{equation}
The nodes' health states are gathered into an $n$-dimensional vector $X(t)\in\{0,1\}^n$, which represents the health state of the population. 

Most of the literature on stochastic epidemic models relies on the assumption that the evolution of the epidemic process $X(t)$ can be represented by a Markov jump process, where the state transitions are triggered by Poisson clocks. Such an assumption, although simplistic, allows to use the rich theory on Markov processes~\cite{levin2006book} to perform rigorous analyses of the model, determining its asymptotic and transient behavior, as detailed in the following. 

\begin{figure}
\centering
\includegraphics[scale=1.8]{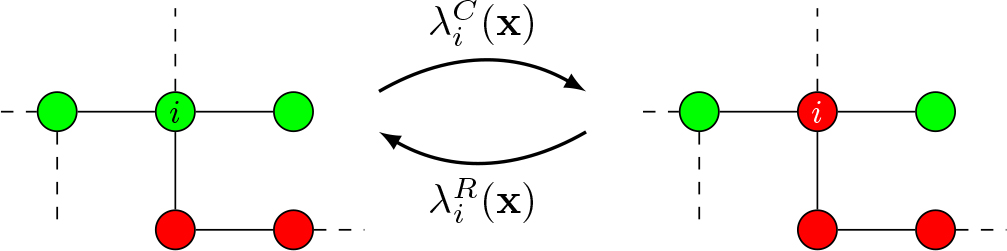}
\caption{Transitions of the Markov process $X(t)$ for a stochastic network SIS model are determined by the contagion and recovery mechanisms. Nodes in green are susceptible and nodes in red are infected.}
\label{fig:markov}\end{figure}

We assume that the vector $X(t)$ evolves according to a continuous-time Markov process, governed by the contagion and the recovery mechanisms, which act on the health state of each node. The former regulates a node's state transitions from the susceptible state to the infected and is modeled by a Poisson clock with the rate proportional to the sum of the connection strengths of the infected individuals that node $i$ is in contact with and to the contagion rate $\lambda_i$, that is,
\begin{equation}
    \lambda^C_i\big(X(t)\big)=\lambda_i\sum_{j\in\mathcal V}a_{ij}X_j(t)\,.
\end{equation}
The recovery mechanism determines the state transitions from the infected state to the susceptible and depends only on the recovery rate of the individual, that is, it is a Poisson clock with the rate equal to $\lambda^R_i\big(X(t)\big)=\mu_i$. These two mechanisms yield the following transition probabilities
\begin{equation}\label{eq:sis_network_stochastic}
    \begin{array}{l}
    \mathbb{P}[X_i(t+\Delta t)=1\,|\,X_i(t)=0]= \lambda^C_i\big(X(t)\big)\Delta T+o(\Delta T),\\
        \mathbb{P}[X_i(t+\Delta t)=0\,|\,X_i(t)=1]= \lambda^R_i\big(X(t)\big)\Delta T+o(\Delta T),
    \end{array}
\end{equation}
for all $i\in\mathcal V$, which unequivocally determine the dynamics of the Markov process, as illustrated in Fig.~\ref{fig:markov}. 

At this stage, one may observe the presence of strong similarities between the expressions of the transition probabilities  of the Markov process and the differential equations that govern the deterministic network SIS model. Indeed, these two models have strong connections. In fact, recently, a different formalization of the deterministic network model has been proposed. In this formalization, each node of the network is a single individual of the population (similar to the stochastic framework) and the corresponding health state $x_i(t)$ denotes the probability that node $i$ is infected at time $t$, that is, $x_i=\mathbb E[X_i(t)]$~\cite{vanmieghem2009,Sahneh2013}. Under the so-called N-intertwined mean-field approach, the following approximation is made: $\mathbb E[X_i(t)X_j(t)]\approx \mathbb E [X_i(t)]\mathbb E[X_j(t)]$, showing that the deterministic system of ODEs in \eqref{eq:sis_network_matrix} approximates the evolution of the expected value of the stochastic process $X(t)$. This approach, even though not exact (since it is based on an approximation that typically does not hold true), is often used to approximate and study the evolution of more complex stochastic network epidemic models.

From \eqref{eq:sis_network_stochastic}, we observe that the disease free equilibrium  $\vect x=\vect 0$ is the unique absorbing state of the Markov process and is globally reachable. Hence, different from its deterministic counterpart, in the stochastic SIS model, the disease is always eradicated with probability $1$ in finite time~\cite{levin2006book}. However, in the following one will see that an epidemic threshold is still present in terms of the duration of the transient evolution of the epidemic outbreak, before reaching the disease-free equilibrium. Formally, it has been observed that a phase transition with respect to the \emph{eradication time} can be established, where the latter is defined by
\begin{equation}
    T:=\min\{t\geq 0:X(t)=\vect 0\}\,.
\end{equation}
Specifically, a set of results that characterize the expected value of the eradication time $\mathbb{E}[T]$ depending on the model parameters and on the network structure has been established~\cite{Ganesh2005,mieghem2013decay,Mountford2016}. The following theorem gathers some key results from the cited literature.
\begin{theorem}\label{th:sis_network_stochastic}
Consider the homogeneous  stochastic network SIS model in \eqref{eq:sis_network_stochastic} with $\lambda_i=\lambda$ and $\mu_i=\mu$, for all $i\in\mathcal V$, on a strongly connected network. If \begin{equation}\label{eq:r0_sto}
    \frac{\lambda}{\mu}<\frac{1}{\rho(A)}\,,
\end{equation}where $\rho(A)$ is the spectral radius of the (weighted) adjacency matrix $A$, then 
\begin{equation}\label{eq:fast_extinction}
    \mathbb{E}[T]\leq\frac{\ln n}{\mu-\lambda\rho(A)}\,;
\end{equation}
if $\lambda/\mu> 1/\rho(A)$, then 
\begin{equation}\label{eq:endemic_disease}
    \mathbb{E}[T]\geq K_1 e^{K_2 n}\,,
\end{equation}
where $K_1,K_2>0$ are two nonnegative constants that depend on the model parameters and on the network structures.
\end{theorem}

The proofs of these results are quite technical and based on the theory of Markov processes~\cite{levin2006book}. Briefly, the key idea is that, in the fast extinction regime, the probability that the number of infected individuals in the population increases as determined by \eqref{eq:sis_network_stochastic} is always less than the probability that it decreases, yielding a drift in the direction of the disease-free equilibrium. Above the epidemic threshold, instead, the inequality is reversed when the process is close to the disease-free equilibrium, implying that the infections tend to rise and large stochastic fluctuations are needed to reach the disease-free equilibrium. Similar techniques have been used to study other stochastic epidemic models on networks, including the SIR model. For more details, we refer to these review papers~\cite{Nowzari2016,Draief2010}.

The results summarized in Theorem \ref{th:sis_network_stochastic} show the existence of a sharp phase transition between a regime where the epidemic is on average quickly eradicated, and a regime where the disease lasts on average for an exponentially long time. However, these results on the average eradication time may fail in characterizing the actual behavior of a single instance of the epidemic process. On the one hand, quick eradication can be guaranteed by directly applying the Markov inequality to \eqref{eq:fast_extinction}, yielding $T\leq (\ln n)^\alpha$ with probability converging to $1$ as $n\to\infty$, for any $\alpha>1$. On the other hand, an accurate analysis of the eradication time has been performed in~\cite{Fagnani2019}, where it has been established that, if $\frac{\lambda}{\mu}>\frac{\sigma}{\rho(A)}$, where $\sigma\geq 1$ is a constant that depends on the network structure (more specifically, on the isoperimetric constants of the network), then the eradication time is exponentially large in the number of nodes with probability converging to $1$ as $n\to\infty$. Note that, if $\sigma>1$, there is a gap between the two regimes. In some specific cases (e.g., for complete graphs and Erdős-Rényi random graphs), $\sigma$ may be equal to $1$, yielding a sharp phase transition not only in the expected duration of the epidemic disease, but also in its actual duration, with high probability.

Recent efforts have been made to relax the Markovianity assumption, assuming different forms for the process $X(t)$. In particular, an SIS model, in which the statistical distribution of the contagion time and/or of the recovery time differ from the exponential distribution associated with Markov processes, was proposed in~\cite{Cator2013}. Therein, a mean-field approach is used to determine conditions for fast eradication of the disease. Further extensions of this approach can be found in~\cite{Kiss2015,Liu2018}. Without relying on any mean-field approximations, in~\cite{Ogura2020}, a lower bound on the decay rate to the infection-free equilibrium is rigorously computed. Equivalences and differences between Markovian and non-Markovian epidemic models have been extensively discussed in~\cite{Starnini2017,Feng2019}.  All these works suggest that the non-Markovianity of the mechanisms that govern the epidemic process may have a large impact on the spread of a disease and outline an important avenue of future research in the field of stochastic epidemic models toward shedding lights on how the distributions of infections and recovery times shape the spreading process.

\begin{table}
\centering
\caption{Notation for network epidemic models}
\begin{tabular}{r|l}
Symbol& Description\\
\hline
$\mathcal G=(\mathcal V,\mathcal E,A)$  & network (node set, node set, adjacency matrix)\\
$\rho(A)$&spectral radius of the adjacency matrix $A$\\
$n$  & number of nodes\\
$X_i(t)$  & health state of node $i$ (stochastic models)\\
$\lambda_i$  & infection rate of node $i$\\
$\mu_i$ & recovery rate of node $i$\\
$\vect x^*$& endemic equilibrium\\
$T$ & eradication time (stochastic models)\\
\end{tabular}

\label{tab:notation}
\end{table}

\subsection{Discrete-time epidemic models}

This survey focuses mostly on continuous-time epidemic models. However, it is important to mention that the continuous-time formulations of deterministic and stochastic models presented in this survey (both the population and the network models) naturally have discrete-time counterparts, where differential equations and Markov processes are replaced by difference equations and Markov chains, respectively. Here, we report the equations for the discrete-time deterministic network SIS model, which will be used in some of the models of epidemics on dynamic networks presented in Section \ref{sec:dynamic}. For each node $i\in\mathcal V$, the health state is updated as follows: 
\begin{equation}\label{eq:sis_discrete}
   x_i(t+1)=(1-\mu_i)x_i(t)+\big(1-x_i(t)\big)\left(1-(1-\lambda_i)^{m_{i}(t)}\right),
\end{equation}
where $m_i(t)$ is defined in \eqref{eq:contagion}; here, $\lambda_i$ and $\mu_i$ have to be interpreted as the per-contact infection probability and the per-time-unit recovery probability, respectively.

Most of the results discussed in this section concerning the existence of a phase transition between a fast extinction regime and a regime where the epidemics becomes endemic and its dependence on the model parameters and on the network structure can be extended with some minor adjustments to the discrete-time counterparts of the models. For more details on the analysis of discrete-time epidemic models and their main results for deterministic models, we refer to~\cite{Allen2008,Brauer2010,Pare2020,liu2020,Prasse2020discrete}, for which a recent review paper by Parè et al. covers most of the results~\cite{Pare2020review}; for stochastic models, we refer to~\cite{Gomez2010,Ahn2014}. Detailed discussions on the main differences between continuous-time epidemic models and their discrete-time counterparts can be found in~\cite{Chan2019,Fennell2016}.

\subsection{Challenges for Network Epidemic Models} 

A few years ago, L. Pellis et al. in a perspective paper outlined eight important challenges for network epidemic models~\cite{Pellis2015}. Besides other --- more practical --- directions, calling for the integration of network computational modeling and epidemiological relevant data, two key challenges were identified, which are of great interest for the engineering community. The first challenge concerns the study of epidemic models on dynamic network structures, leading to the study of a nonlinear time-varying dynamical system. The second focuses on understanding how the network structure (static or dynamic) can be exploited to effectively design intervention policies to stop or mitigate the disease spreading; control-theoretic tools are key to address this second challenge. In the rest of this survey, we focus on these two research directions, presenting the state of the art in terms of key progresses of the last few years and most promising lines of current research.

\section{Epidemics models on dynamic networks}\label{sec:dynamic}

The extension of the classic compartmental models to static networks and the corresponding rigorous analysis have allowed the scientific community to understand how the architecture of human social interactions affects the spread of epidemic diseases in interconnected populations. However, in most of the real-world epidemic outbreaks, the underlying network of social interactions is not static, but dynamically changes, co-evolving with and influenced by the spread of the disease~\cite{Gross2008,Bansal2010,Holme2012,Holme2015}. 

Several endogenous and exogenous reasons may be adduced to explain and motivate the dynamic evolution of the network structure. First, dynamical changes of the individuals' patterns of interactions may be directly or indirectly caused by seasonal factors, such as school holidays and weather conditions, which may favor or hinder gatherings and social events~\cite{Altizer2016,Bansal2010}. Second, social interactions are often characterized by an intermittent behavior, whereby individuals' propensity to generate connections is subject to burstiness, yielding clusters of connections separated by latency periods~\cite{barabasi2005, Goh_2008}. Third, the infection events themselves may affect the network structure. In fact, not only infected individuals may reduce their interactions as a consequence of the illness, but also susceptible individuals --- driven by awareness and risk perception cognitive mechanisms --- may dynamically modify their behavior, reducing or rewiring their interactions to avoid contagion~\cite{Gross2008,Funk2010, Verelst2016}. Finally, nonpharmaceutical intervention policies often entail a dynamical modification of the pattern of human interactions, which may be dynamically reduced through the implementation of lockdown or social distancing policies, or reshaped by travel bans and mobility limitation, as observed during the ongoing COVID-19 pandemic~\cite{Lai2020}. All these evidences of the dynamic nature of social interactions call for the extension of the epidemic models presented in the previous sections to dynamic networks. 

Here, we review and discuss some successful modeling paradigms that have been recently developed to capture this dynamic nature of the network of social interactions within analytically tractable mathematical models of epidemic diseases. Through the analysis of these models, we shed lights on how the dynamical nature of human interactions plays a key role in shaping the evolution of the epidemic outbreak. As we shall illustrate in the next section, the insight gained into the epidemic process through these analyses has enabled researchers  to propose and assess valuable intervention strategies to control and mitigate the epidemic spreading, taking into account and leveraging the dynamical properties of social interactions.

\subsection{First approaches: time-scale separation}

\begin{figure}
\centering
\includegraphics[scale=1.8]{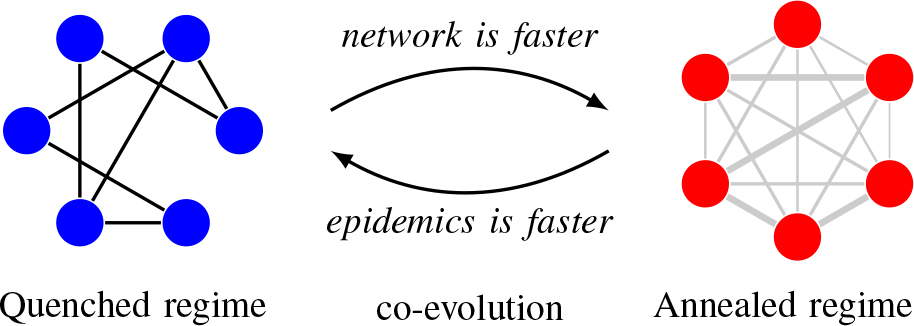}
\caption{Modeling paradigms for dynamic networks. Under the assumption of time-scale separation, quenched and annealed regimes can be found. Between these regimes, several paradigms have been proposed to capture the co-evolution of the epidemic process and the network structure at comparable time-scales, including temporal-switching, activity-driven, and edge-Markovian networks, detailed in this survey. }
\label{fig:temporal}
\end{figure}

The first class of approaches to deal with dynamic networks has extensively relied on time-scale separation techniques. These techniques are based on the assumption that the epidemic process and the network dynamics evolve at different paces, as illustrated in Fig.~\ref{fig:temporal}. If the epidemic process evolves much faster than the network of interactions, the system is in the so called \emph{quenched regime}. In this regime, static networks are accurate proxies of slowly switching topologies, and the corresponding results presented and discussed in the previous section are thus used to study the evolution of the epidemic outbreak. On the other extreme, we encounter the \emph{annealed regime}, in which the evolution of the network is assumed to be much faster than the disease spreading process~\cite{Newman2002}. In this regime, if the following limit exists
\begin{equation}\label{eq:barA}
 \bar{A}=\left\{\begin{array}{ll}\displaystyle\lim_{T\to\infty}\frac{1}{T}\sum_{t=1}^T A(t)&\text{ if }t\in\mathbb Z^+\,,\\
 \displaystyle\lim_{T\to\infty}\frac{1}{T}\int_{t=0}^T A(t)dt&\text{ if }t\in\mathbb R^+\,,
 \end{array}\right.
\end{equation}then we can define the average graph $\mathcal G=(\mathcal V,\bar{\mathcal E}, \bar A)$,  with $(i,j)\in\bar{\mathcal E}\iff \bar a_{ij}>0$, and the dynamic network can be effectively represented and studied by means of its average graph. Also in the annealed regime, results on static networks are applied to the average graph to study the behavior of the dynamical system.

In the physics community, epidemics on annealed networks have been extensively studied, aiming at computing --- or approximating --- the epidemic threshold of the model on a network with known degree distribution, without going through the explicit computation of the spectral radius of the average adjacency matrix. Specifically, by relying on a degree-based mean-field approach, the epidemic threshold for the SIS model on unweighted networks has been computed in~\cite{Pastor-Satorras2001} as a function of the degree distribution of the network. Note that such an approximation is exact in the limit of large-scale systems $n\to\infty$. Specifically, let 
\begin{equation}\langle \bar k\rangle=\frac{1}{n}\sum_{i\in\mathcal V}\bar k_i\,,\qquad \langle \bar k^2\rangle=\frac{1}{n}\sum_{i\in\mathcal V}\bar k_i^2\,,\end{equation} 
denote the mean and the second moment of the degree distribution of the average network $\bar{\mathcal G}=(\mathcal V, \bar{\mathcal E},\bar A)$, respectively. Then, the following expression for the epidemic threshold can be obtained for uncorrelated annealed networks (that is, if $\bar a_{ij}\propto \bar k_i\bar k_j$). 
\begin{theorem}\label{th:sis_annealed}
Consider an SIS model on a dynamic network in the uncorrelated annealed regime with $\lambda_i=\lambda$ and $\mu_i=\mu$, for all $i\in\mathcal V$. Let us define the following epidemic threshold:
\begin{equation}\label{eq:sis_annealed}
\sigma=\frac{\langle \bar k\rangle}{ \langle \bar k^2\rangle}\,.
\end{equation}
Then, in the limit $n\to\infty$, if ${\lambda}/{\mu}<\sigma$, the disease-free equilibrium is asymptotically stable; otherwise, if ${\lambda}/{\mu}>\sigma$, the disease-free equilibrium is unstable. 
\end{theorem}
From Theorem \ref{th:sis_annealed}, we establish that, below the epidemic threshold in \eqref{eq:sis_annealed}, the epidemic is in the fast extinction regime while, above the epidemic threshold, the epidemic becomes endemic. An important implication of \eqref{eq:sis_annealed} applies to scale-free networks, whose degrees follow a power-law distribution. In fact, in many real-world applications, complex networked systems are modeled by scale-free networks with the power-law exponent  between $2$  and $3$~\cite{Barabasi1999}. In these scenarios, the expression of $\sigma$ in \eqref{eq:sis_annealed} vanishes as $n\to\infty$, implying that for any nonzero contagion rate $\lambda>0$, epidemics always spread on large scale-free networks. More details on this approach and on further extensions of these works to more general networks (including correlated networks) and to more complex epidemic models (including the SIR model) can be found in~\cite{RevModPhys.87.925}.

The quenched and annealed regimes discussed in the above rely on the assumption that the epidemic process and the network evolve at different paces and, thus, on different time-scales. However, the arguments raised at the beginning of this section to motivate the need for dynamic networks provide evidence that such a time-scale separation is often restrictive and unrealistic, since the contagion process and the network evolution are often intertwined and thus often evolve at comparable time-scales. In the last few years, several efforts have been made to overcome the limitation of time-scale separation and propose a theory for epidemics on dynamic networks which allow to model and study the coevolution of the two dynamical processes. In the rest of this section, we present and discuss some of these fascinating paradigms.

\subsection{Temporal-switching networks} 

The use of temporal-switching networks to study epidemics in time-varying systems has been initially proposed in~\cite{Prakash2010}. In its original incarnation, a switching network is generated by repeating a deterministic sequence of $T$ static networks, characterized by the adjacency matrices $A_1,\dots,A_T$, so that $A(t)=A_{t \text{ mod } T}$. Therein, the discrete-time deterministic SIS model has been studied, extending the results found for static networks. Specifically (in the homogeneous scenario), if one defines 
\begin{equation}
    P:=\prod\limits_{t=1}^T\Big((1-\mu)I+\lambda A_t\Big)\,,
\end{equation}
then, the behavior of the system is determined by the spectral radius $\rho(P)$. If $\rho(P)<1$, then the disease-free equilibrium is asymptotically stable and the disease is quickly eradicated. Otherwise, the disease-free equilibrium becomes unstable and the disease becomes endemic. Such a framework has been extended to the discrete-time stochastic SIS model in~\cite{Valdano2015}, showing  the same epidemic threshold by mapping the time-varying system onto a multi-layer network structure. In~\cite{Sanatkar2016}, sufficient conditions for stability have been established for a more general scenario, where the networks switch arbitrarily among a set of topologies, possibly according to  stochastic mechanisms, such as Markov switching rules where, given a Markov process $\sigma(t)$ with the state space $A_1,\dots,A_T$, we set $A(t)=A_{\sigma(t)}$. This sufficient condition is expressed in terms of the maximum
possible norm of products of matrices $P_t$ in the set $$\mathcal P:=\left\{P_t=(1-\mu)I+\lambda A_t\right\}\,.$$
For a review of the most recent developments of this theory, including the computation of a unified formula for the epidemic threshold, we refer to the following paper by Zhang et al.~\cite{Zhang2020}.

The use of temporal-switching networks to study epidemics on dynamic networks has been recently extended to continuous-time processes. Specifically, in~\cite{Rami2014}, the theory of positive linear switched systems is leveraged to derive conditions for global asymptotic stability of the disease-free equilibrium. Such conditions are obtained by combining the stability analysis of the linearized switched system an appropriate notion of irreducibility for the linearization. Specific results are obtained if the topology evolves stocastically according to a Markov switching rule. This approach is followed in~\cite{Ogura2017switching} to design Markov switching laws to enforce quick eradication of the disease via geometric programming. In~\cite{Speidel2016}, a continuous-time network SIS model is studied in a scenario where the network topology switches deterministically, at discrete-time steps, following a sequence of adjacency matrices. Therein, the epidemic threshold is computed in terms of the Lie commutator bracket of the adjacency matrices, showing that adjacency matrices that are non-commuting yield a lower epidemic threshold, favoring the epidemic spreading. In~\cite{Valdano2018}, the scenario of continuous-time switching networks is considered. The epidemic threshold is explicitly computed in the case when the adjacency matrix $A(t)$ commutes with the aggregated adjacency matrix up to that time $\int_0^tA(s)ds$. Under this condition, the order of the switching matrix has no effect on the dynamics, which is fully determined by the average adjacency matrix $\bar A$, defined in \eqref{eq:barA}. If ${\lambda}/{\mu}<1/\rho(\bar A)$, then the disease-free equilibrium is asymptotically stable; otherwise, it is unstable and the disease becomes endemic. In~\cite{Pare2018} the analysis of the continuous-time deterministic heterogeneous SIS model from~\cite{Wang2003} is extended to slowly changing time–varying systems by leveraging Lyapunov arguments. In that work, the effect of stochastic perturbations of the model is also studied.

A major limitation of the theoretical analysis of temporal-switching networks is that it is often assumed that the network switches between different adjacency matrices determined a-priori and that the switches typically take place at deterministic (often equally-spaced) time-instants. In the following, we will present two alternative paradigms that do not rely on this assumption. 

\subsection{Activity-driven networks}\label{sec:adn}

\begin{figure}
\begin{center}
\subfloat[$t=0$]{\includegraphics[scale=1.8]{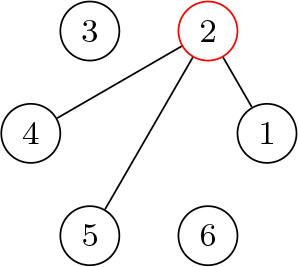}}\qquad
\subfloat[$t=1$]{\includegraphics[scale=1.8]{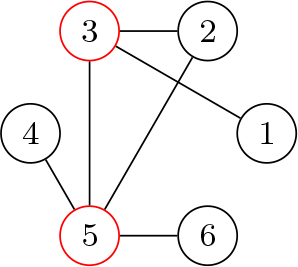}}\qquad
\subfloat[$t=2$]{\includegraphics[scale=1.8]{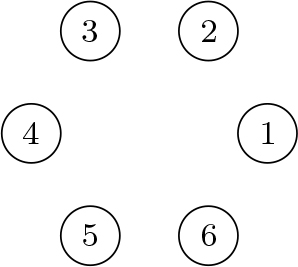}}
\end{center}
\caption{Exemplary evolution of a discrete-time ADN. At discrete-time $t=0$, node $2$ is activated  and generates $m=3$ undirected links. At discrete-time $t=1$, nodes $3$ and $5$ {are activated}, generating $m=3$ undirected links each. At discrete-time $t=2$ none of the nodes is activated and, consequently, no links are generated.}
\label{fig:adn}
\end{figure}

In all the network models presented so far, interactions were determined by some pre-determined connectivity patterns, represented by an adjacency matrix, or a sequence of adjacency matrices. Activity-driven networks (ADNs), originally proposed by N. Perra et al. in~\cite{Perra2012} suggested to perform a paradigm shift. In ADNs, interactions are not seen as a consequence of a network structure which identifies pre-determined dyadic relationships between pair of nodes, but are rather generated by individuals' properties, according to a stochastic process. In their original incarnation, each individual $i$ is characterized by a unique parameter $a_i\in[0,1]$, termed \emph{activity}, which represents an individual's propensity to generate interactions. Then, starting from $t=0$, the dynamic network is generated according to the following algorithm:
\begin{enumerate}[i)]
    \item the edge set is initialized as $\mathcal E(t)=\emptyset$;
    \item each individual $i$ activates  with probability $a_i$, independent of the others. If an individual $i$ activates, then the individual generates $m$ links with $m$ other individuals $\{i_1,\dots, i_m\}$, chosen uniformly at random as an $m$-tuple in the population;
    \item undirected links $(i,i_h)$, $h=1,\dots, m$ are added to the edge set $\mathcal E(t)$; and
    \item the time index is increased by 1 step, and the algorithm resumes from step i).
\end{enumerate}
The formation process of an ADN is depicted in Fig.~\ref{fig:adn}. Note that, at each discrete-time, the network generated by the ADN process is always undirected, but not necessarily connected.

The main strength of this paradigm lies in its simplicity: the temporal nature of the system is encapsulated in a single $n$-dimensional vector. Such a simple formulation has enabled to perform rigorous analytical studies of the properties of the network generated and of dynamical processes evolving on it. Specifically, in~\cite{Perra2012}, the following result has been established.
\begin{theorem}\label{th:sis_adn}
Consider an SIS model on an ADN with $\lambda_i=\lambda$ and $\mu_i=\mu$, for all $i\in\mathcal V$. Let us define the following epidemic threshold:
\begin{equation}\label{eq:sis_adn}
    \sigma=\frac{1}{m\big(\langle a\rangle+\sqrt{\langle a^2\rangle}\big)}\,.
\end{equation}
For $\lambda/\mu<\sigma$ the epidemic is quickly eradicated with probability converging to $1$ as the network size $n\to\infty$, while for $\lambda/\mu>\sigma$ the epidemic becomes endemic with probability converging to $1$ as the network size $n\to\infty$.
\end{theorem}

The original formulation of ADNs was proposed in a discrete-time framework, so the discrete-time counterpart of the SIS model in \eqref{eq:sis_discrete} was studied. A continuous-time formulation of ADNs has been proposed in~\cite{Zino2016}, where the synchronous activations ruled by the activity-based mechanism are substituted by an asynchronous mechanism, where each node is activated according to a Poisson process with the rate equal to its activity, yielding thus a Markov process. Therein, the continuous-time SIS model is analyzed, giving rise to the same threshold identified in Theorem \ref{th:sis_adn}.

One of the major advantages of the ADN formulation is its simplicity that, besides enabling rigorous analytical studies, allows to expand the formalism in several directions. In fact, in the last few years, several extensions and generalizations of the ADN paradigm have been proposed. These extensions allow to include several features of real-world complex systems in the model. The analytical tractability of ADN-based models has enabled the rigorous computation of the epidemic threshold for these models and the characterization of their behavior, similar to Theorem \ref{th:sis_adn}, shedding light on the role of these real-world phenomena on the spreading of epidemics. These extensions include the presence of preferential connectivity patterns~\cite{sun2015contrasting,Lei2016, Nadini2020}, community structures~\cite{Nadini2018,Bongiorno2019}, heterogeneous propensity to receive connections~\cite{Pozzana2017}, memory and burstiness in the link formation process~\cite{PhysRevLett.114.108701,Zino2018,Mancastroppa2019}, and high-order relations~\cite{Petri2018}. A detailed review of the results for classical activity-driven networks and for the explicit results of the epidemic thresholds for these recent extensions can be found in~\cite{Leitch2019}.

\subsection{Edge-Markovian dynamic graphs}

\begin{figure}
\begin{center}
\includegraphics[scale=1.8]{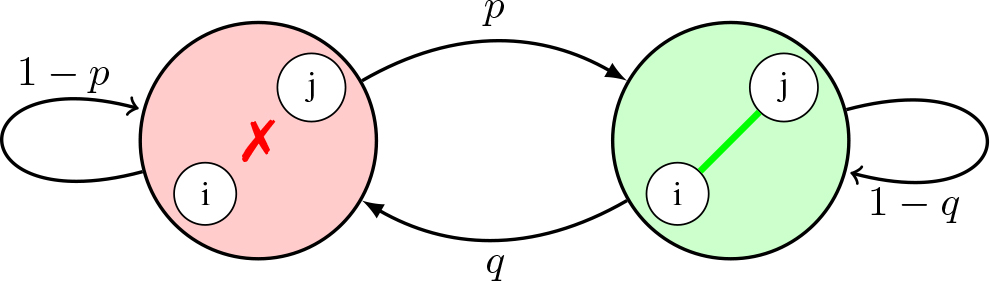}
\end{center}
\caption{Transitions of the two-state Markov chain associated with the existence of a generic link $(i,j)$ of a discrete-time edge-Markovian dynamic graph.}
\label{fig:edge_markov}
\end{figure}

A different approach, which to a certain extent combines the presence of a connectivity pattern, determined a-priori, and the stochasticity of its evolution, are edge-Markovian dynamic graphs. This paradigm has been proposed in~\cite{Clementi2010} to model stochastic evolution of dynamic networks. In edge-Markovian dynamic graphs, each potential link (edge) of the graph (i.e., each pair $(i,j)\in\mathcal V\times\mathcal V$, with $i\neq j$) is associated with a two-state Markov chain (independent of the other links), where the two states represent the existence and nonexistence of the link, respectively. Two probabilities $p,q\in[0,1]$ are defined so that, at each discrete time-step, the chain switches from nonexistence to existence with probability $p$, while the opposite transition happens with probability $q$, as illustrated in Fig.~\ref{fig:edge_markov}. In plain words, the network is initialized as a given (static) network. Then, each link that exists at time $t$ disappears at the following time-step with probability $q$ (independent of the others), while nonexisting links at time $t$ appear with probability $p$ (independent of the others). A continuous-time formulation of the model can be obtained by substituting the Markov chains with continuous-time Markov processes~\cite{Kiss2012}.

Epidemic processes on edge-Markovian dynamic graphs have been proposed and studied in~\cite{Kiss2012, Taylor2012}. Specifically, in~\cite{Taylor2012}, the authors derived the value of the epidemic threshold in terms of the basic reproduction number, which has a complex expression. Edge-Markovian dynamic graphs are amenable of several analytically tractable extensions to overcome the limitation posed by its original formulation. For instance, the Markov chain (or process) underlying the network evolution determines the inter-event time distribution for the appearance and disappearance of the links (geometric in the discrete-time model, exponential in the continuous-time model), which is not consistent with the presence of burstiness and temporal clustering discussed in the above~\cite{barabasi2005,Goh_2008}. In~\cite{Ogura2016}, this paradigm has been extended to account for general inter-event time distributions, establishing a sufficient condition for (almost sure) exponential stability of the disease-free equilibrium.

\section{Control of epidemics on networks}\label{sec:control}

For engineering researchers who are familiar with the monitoring, intervention or control of complex systems, it is of great interest to know what research works have been carried out to study how to influence, mitigate, and even stop the epidemic processes, especially those based on the models we have presented in the previous sections. Although much fewer control results have been produced compared to the epidemic modeling activities, it is beneficial to give an overview of the existing results, which will serve to inspire researchers, with or without control theory background, to work in this critically important research area. In what follows, we categorize the corresponding results into control of deterministic and stochastic epidemics respectively, and to underscore the importance, devote a separate subsection for the discussion of the distinct features when the underlying networks are dynamic. Note that an earlier survey~\cite{Nowzari2016} has summarized some main results on control epidemics on networks by then, and thus we give special attention, wherever appropriate, to those that appeared in the last five years.

\subsection{Control of deterministic epidemics on static networks}

\begin{figure}
\begin{center}
\includegraphics[scale=1.8]{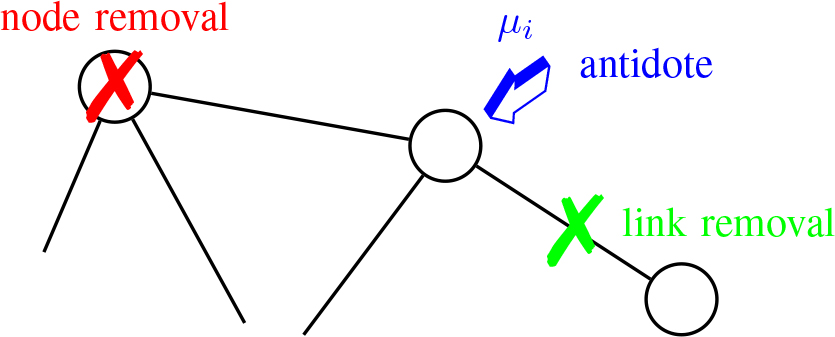}
\end{center}
\caption{Schematic of the main control actions that can be taken in the control of static networks.}
\label{fig:control}
\end{figure}

The simplest idea of controlling epidemics processes on networks comes from the intuition that removing infected or high-risk individuals and the links associated with them will slow down the propagation, which in practice translates to quarantine and vaccination policies. Following the discussion in the previous sections, this intuition implies that one can lower the epidemic threshold, e.g., by reducing the spectral radius $\rho(A)$ of the adjacency matrix $A$. Another intuitive idea is to optimize the distribution of antidote, which in practice translates into modifying the entries of matrix $M$ in \eqref{eq:sis_network_matrix}. These key control actions are illustrated in Fig.~\ref{fig:control}. Naturally, there is always a cost associated with the control actions and consequently optimal or sub-optimal control objectives can be formulated. However, to get a flavor on why such intuitive control problems under cost constraints are difficult to solve, we look at the following direct formulation of the control problem. 

\begin{problem}\label{problem:removal}Given a network of $|\mathcal E|=\ell$ links and a fixed number $m<\ell$. Find those $m$ links $\mathcal R=\{e_1,\dots, e_m\}$ of the network after removing which the adjacency matrix of the resulting network $\tilde G=(\mathcal V,\mathcal E\setminus \mathcal R)$ has the minimum spectral radius among all the reduced networks obtained by removing $m$ links from the original network.
\end{problem}

Although this control problem is straightforward to formulate, it is very difficult to be solved analytically. In fact, it was shown in~\cite{vanmieghem2011} that Problem \ref{problem:removal} is NP-hard. A similar spectral minimization control problem through removing nodes is discussed in~\cite{Preciado2009}. In fact, the difficulty in solving such control problems is rooted in the fact that the formulated optimal control problems are variations of the constrained combinatorial optimization problems, which are in general hard problems. For this reason, various heuristics have been proposed to solve the control problems, and a lot of them have smartly taken advantage of the network structures, e.g., the degree distribution of the nodes, centrality indices, and connectivity patterns~\cite{PhysRevE.65.036104,chung2009}, or solving a non-convex quadratically constrained quadratic program~\cite{enns2012}.

Another intuitive approach lies in tuning the values of the parameters by increasing the recovery rates $\mu_i$  to minimize the expression in \eqref{eq:threshold_het}, which characterizes the epidemic threshold. Such an approach can be interpreted as a resource allocation problem, where limited amounts of antidote can be distributed to the population. Such an intuition is formalized in the following optimization problem.
\begin{problem}\label{problem:resource}Given a network and a fixed budget $B>0$
\begin{equation}
    \begin{array}{rl}
       \text{minimize}_{\vect\mu}&\rho(\Lambda A-M)\\
         \text{subject to}& f(\vect\mu)\leq B\,,\\
         & \mu_i\in[\underline \mu_i,\overline\mu_i],\,\,\forall\,i\in\mathcal V\,,
    \end{array}
\end{equation}
where $f:\mathbb{R}^n\to \mathbb{R}^+$ is a cost function associated with the cost of increasing the recovery rate and $\underline \mu_i$ ($\overline \mu_i$) is the minimum (maximum) admissible recovery rate for node $i\in\mathcal V$.
\end{problem}
Even though the objective function is a spectral radius, which in general is non-convex, under reasonable assumptions on the cost function, tools from geometric programming and convex optimization can be leveraged to tackle the problem. Solutions have been proposed in a centralized fashion~\cite{Gourdin2011,Preciado2014,Nowzari2015,Ottaviano2017}, and through distributed approaches~\cite{Enyioha2015distributed,ramirez2018,Mai2018,Mai2019,Somers2021}. Some of these works deal with a more general problem, in which, besides increasing the recovery rate, the controller can also reduce the infection rates $\lambda_i$ modeling, for instance, the distribution of personal protective equipment. A resource allocation problem similar to Problem \ref{problem:resource} is studied for an extension of the SIS model, in which complication phenomena of the illness are considered~\cite{Hung2019,DiGiamberardino2019}.

When considering more complex dynamics than the standard SIS epidemic model, ideas from  optimal control have already been applied a decade ago~\cite{Bloem2009}, using a linear-quadratic regulator, and, more recently, in~\cite{Eshghi2016}, leveraging the Pontryagin’s Maximum Principle.

More recently, researchers have identified impossibility results which reveal the possible limitation of feedback control. In~\cite{Ye2020,Ye2021}, the authors have proved that, utilizing the recovery rate $\mu_i$ as control input, a large class of  distributed  controllers cannot guarantee convergence to the disease-free equilibrium. In~\cite{Liu2019}, a similar result has been proved for more complex dynamics that involve two concurrent epidemic processes. The limitations may become even more profound when examining the effect of optimal control in real-world disease management~\cite{Bussell2019}.

\subsection{Control of stochastic models on static networks}

Because of the stochastic nature of the models, the related control results are centered around evaluating the control performance in terms of bounding as tightly as possible the epidemic thresholds on different classes of networks.  In~\cite{Borgs2010,Drakopoulos2014}, the studied control problem is how to distribute a fixed amount of antidote to nodes in the given network that may have special topological features, e.g., scale-free networks. In~\cite{Borgs2010},  two different methods are compared, one based on contact tracing, which augments the recovery rates of all neighbors of an infected node, and the other based on degree-centrality,  which augments the recovery rates of all nodes, proportional to their degrees. Surprisingly, it is found that contact tracing may only succeed when the number of infected individuals is small (e.g., in early stages of the epidemic outbreak), since it requires a total amount of antidote $B=\sum_{i}\mu_i$ that grows super-linearly in the number of contacts; otherwise the degree-centrality based approach outperforms contact tracing, as stated in the following result.
\begin{theorem}\label{th:borgs}
Consider the stochastic network SIS model in \eqref{eq:sis_network_stochastic} on a generic network. If $\mu_i\geq\lambda_ik_i$, then the expected eradication time verifies $\mathbb{E}[T]\leq K\ln n$, for some constant $K>0$.
\end{theorem}
However, from this result, we note that the needed amount of antidote $B$ scales linearly with the sum of the weights of the links in the network, hindering its practical implementation even in sparse networks, where such a sum grows linearly in the number of nodes.

In a similar setup, in~\cite{Drakopoulos2014, Drakopoulos2017}, another control method is proposed, in which the antidote is dynamically allocated to the nodes. The allocation method in~\cite{Drakopoulos2014} requires that all the antidote $B$ is concentrated at a single node at each time; using martingale theory, the following result is established.
\begin{theorem}
Consider the stochastic network SIS model in \eqref{eq:sis_network_stochastic} on a generic network with the control policy proposed in~\cite{Drakopoulos2014}. Let us define the maximum degree and the cutwidth of the network as \begin{equation}
    \Delta=\max_{i\in\mathcal V} k_i\quad\text{and}\quad W=\min_{\emptyset\subset\mathcal S\subset \mathcal V}\sum_{i\in \mathcal S,j\notin \mathcal S}a_{ij},
\end{equation}
respectively. If $B\geq16\Delta\ln n$ and $B>4W$, then the expected eradication time verifies $\mathbb{E}[T]\leq26n/B$. Moreover, if $B>Kn/\ln n$, then $\mathbb{E}[T]\leq K'\ln n$, for some constants $K,K'>0$.
\end{theorem}
Briefly, the proposed control technique guarantees fast eradication with a sub-linear amount of antidote, depending on the network topology. In~\cite{Drakopoulos2017}, a fundamental limitation is further established for any dynamic allocation by showing how the network topology influences the possibility of eradication the epidemics under the given budget of the antidote to be allocated. A similar control policy based on dynamical resource allocation has been proposed in~\cite{Scaman2016}.

Very recently, the powerful tool of model predictive control (MPC), already used for deterministic epidemic models~\cite{selley2015}, has been used to deal with the control of stochastic epidemic processes without relying on mean-field approximation~\cite{Kohler2018}. The use of MPC has also allowed to deal with more complex epidemic models. For instance, in~\cite{Watkins2020}, the authors deal with a model with a pre-symptomatic phase, by utilizing MPC with a robust moment closure technique. Optimal control has also been considered for stochastic SIS model with different assumptions on the information available~\cite{Grandits2019}.
 
\subsection{Differences in control when the networks are dynamic}

When the networks are dynamic, the optimization problem in control may have to face time-varying constraints. Fortunately, some of the geometric program techniques still apply, although the complexity in seeking the solution increases~\cite{Nowzari2015b}.  Optimal control theory can be applied to time-varying systems; however, it is well known that the corresponding stability conditions might be more conservative and more difficult to check~\cite{Youssef2013}. In~\cite{Ogura2017switching}, the authors deal with a resource allocation problem similar to Problem \ref{problem:resource}, on temporal-switching networks, under the assumption that the switching is determined by a Markov process. Therein, the problem is solved via geometric programming, finding a solution with a computational complexity that grows super-linearly --- but polynomially --- in the network size.

For some specific cases, control strategies for eradicating the outbreak have been successfully proposed. For instance, in~\cite{Gracy2020}, a distributed control scheme has been designed for a deterministic SIS model on periodic time-varying networks. In this scheme, that recalls the antidote distribution proposed in~\cite{Borgs2010} and summarized in Theorem \ref{th:borgs}, each node dynamically sets its recovery rate proportional to the sum of the weight of the links to its neighbors at that time, that is, \begin{equation}
  \mu_i(t)=\lambda_i\sum_{j\in\mathcal V}a_{ij}(t).  
\end{equation}
Under this scheme, global asymptotic stability of the disease-free equilibrium is guaranteed, under the assumption that all the instances of the time-varying network are strongly connected~\cite{Gracy2020}.

An interesting research line is to incorporate the human behavioral mechanisms that lead to variations of the dynamic network structures~\cite{Verelst2016}. Understanding human behavior will be critical to draft and implement vaccination policies~\cite{Liu2014}, which affects the dynamics of the networks and at the same time is deeply affected by the dynamics of the networks~\cite{Liu2014}. Further study may dive in how isolation of infected individuals may adaptively change and reshape the network dynamics and how this can be leveraged to devise effective intervention policies~\cite{RizzoPRE2014,Zino2020ejc}. In particular, in~\cite{RizzoPRE2014}, the authors explicitly compute the following epidemic threshold for a network SIS model on activity-driven networks, depending on the possibility of isolating infectious individuals by decreasing their social activity to a factor $p\in[0,1]$, where $p=1$ means that interventions are enacted and $p=0$ models a complete isolation of infected individuals:
\begin{equation}\label{eq:sis_adn_behavior}
    \sigma=\frac{2}{m\Big((1+p)\langle a\rangle+\sqrt{(1-p)^2\langle a\rangle^2+4p\langle a^2\rangle}\Big)}\,.
\end{equation}
By comparing \eqref{eq:sis_adn_behavior} with \eqref{eq:sis_adn}, one can assess the effect of isolation policies on increasing the epidemic threshold.

More recently, awareness-based control strategies, which were developed a decade ago~\cite{Funk2009}, have been extended to study  temporal networks such as activity-driven networks~\cite{Ogura2019, Yang2019, Zino2020}. The proposed formalism allows to study scenarios in which aware individuals reduce the risk during their physical interactions~\cite{Yang2019} or they dynamically rewire themselves to avoid interactions with infected individuals~\cite{Ogura2019}, possibly in combination with other control policies, such as isolation of infected individuals~\cite{Ogura2019,Zino2020}. In the model proposed by Yang et al. in~\cite{Yang2019}, awareness is modeled as a process that co-evolves with the epidemic spreading on a two-layer network. The epidemic process spreads on a physical contact layer, while awareness spreads on a virtual communication layer. Aware individuals reduce their infection probability, as a consequence of the adoption of self-protective practices. The epidemic threshold is then computed as a function of the individuals' activities on the two layers and the coupling between them. In the activity-driven adaptive-SIS model proposed in~\cite{Ogura2019}, two modifications are made to the ADN algorithm illustrated in Section \ref{sec:adn}. First, in step ii), the activity of an infected individual $i$, is multiplied by a parameter $p_i\in[0,1]$, similar to~\cite{RizzoPRE2014}. Second, in step iii), a link $(i,i_h)$ to a node $i_h$ that is infected is added to the edge set with probability $\pi_i\in[0,1]$. The epidemic threshold is computed, in terms of the decay ratio to the disease-free equilibrium, finding a rather complicated expression that depends on the joint distributions of the activities and the parameters $p_i$ and $\pi_i$. Such an analytic expression is used to devise an optimization problem, to optimally allocate resources into isolation of infected individuals and awareness, solved via geometric programming.

\section{The current challenges with COVID-19}\label{sec:covid}

Since its inception in December 2019, the COVID-19 outbreak has rapidly spread, becoming a worldwide pandemic with over 108 million reported cases as of February 16, 2021. In response to this unprecedented health crisis, we witnessed an extraordinary mobilization of the scientific community toward better understanding the novel disease and combating its spread. Within this joint effort, the systems and control communities are playing a key role in developing accurate mathematical models to predict the evolution of the pandemic and assess the effect of diverse intervention policies that have been enacted or that may be implemented~\cite{Bertozzi2020,Vespignani2020model,SpectrumCovid}.

\begin{figure}
\begin{center}
\includegraphics[scale=1.8]{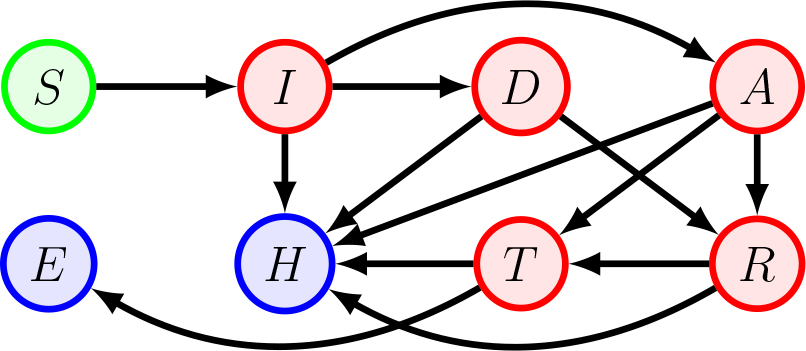}
\end{center}
\caption{Schematic of the SIDARTHE model, proposed in~\cite{Giordano2020} to capture the epidemic progression of COVID-19. The modeled is obtained from an SIR model, in which the health state ``infected" is substituted by five states, representing all the four possible combinations of whether the individual has symptoms and is detected, and a state for individuals seriously ill; two states are used instead of the ``removed" state to account for individuals that are healed ($H$) or extinct ($E$), respectively. The arrows illustrate all the possible transitions, which are governed by $16$ different parameters.}
\label{fig:covid}
\end{figure}

A first, important contribution comes from the formulation and analysis of more complex models for the epidemic progression in fully-mixed populations. These models allow to capture the specific features of COVID-19, including a latency period before the symptoms onset, the presence of asymptomatic individuals, and the implementation of intervention policies such as hospitalization of severe cases and testing. Among these models, we should mention the SIDARTHE model proposed in~\cite{Giordano2020} (see Fig.~\ref{fig:covid}), which takes into consideration also the imperfect reporting of infected individuals. In this work, the epidemic model is studied in a fully-mixed population framework by seeing it as a positive linear system under feedback, and the stability of the disease-free equilibrium is thus characterized in terms of the $H_\infty$ norm of its transfer function, which interestingly coincides with the basic reproduction number $\mathcal R_0$. Using this model, an open-loop fast switching strategy to control and suppress the spread of the disease, consisting in intermittent lockdown phases, has been proposed and studied in~\cite{Bin2020}.

Once the epidemic model has been tailored to capture the epidemic progression of COVID-19, its implementation on a network structure (extensively discussed in this survey) is key to predict the spatial spread of the disease, as highlighted in~\cite{Bertozzi2020}. The time-varying nature of human mobility as well as the implementation of intervention policies through different sequential phases have put epidemic models on dynamic network at the forefront of the stage. In this vein, we mention several data-informed analyses of the outbreak in different countries, using models with regional granularity. These studies include Italy~\cite{Gatto2020,Parino2020}, Ontario, Canada~\cite{Karatayev2020}, Western Australia~\cite{Small2020}, and Kazakhstan~\cite{Kuzdeuov2020}. Based on these network models, non-linear MPC has been adopted to understand the impact of intervention policies and plan their optimal implementation. We mention the works in~\cite{Carli2020} and in~\cite{kohler2020robust}, with case studies based on the outbreak in Italy and Germany, respectively. The impact of local, time-varying lockdown measures and mobility restrictions is analyzed in~\cite{dellarossa2020}, where feedback control laws coordinated by a centralized controller are used to design an intermittent and differentiated regional intervention scheme that outperforms nationwide measures. Dynamic networks also enable modeling the individual behavioral response to the pandemic in terms of the adoption of self-protective behaviors and social-distancing measures, which has a huge impact on disease spreading~\cite{Bavel2020}.

At the moment of writing this survey, effective pharmaceutical treatments for COVID-19 were unfortunately still not available, while the research for a vaccine has recently led to some promising findings, and extensive vaccination campaigns are getting started. In this delicate phase, in which only a very limited amount of the vaccine is available, public health authorities need to carefully plan the distribution of vaccines. The control problems described and discussed in this survey (e.g., how to optimally modify the recovery rates by distributing a fixed amount of antidote) will be precious tools to help inform vaccination campaigns and distribution of pharmaceutical treatments.

\section{Directions of current and future research}\label{sec:future}

In this survey, we went through $260$ years of progresses in mathematical models of epidemics. Starting from the very first intuition that mathematics, and in particular systems theory, can provide effective tools toward better understanding the spread of infectious diseases, we reviewed the recent developments and the state of the art in the analysis and control of epidemic models on networks. Alongside, we outlined some of the most promising avenues of current and future research, which are of particular interest for engineering researchers. These directions include the analysis of epidemics on dynamic networks, where the process is modeled as a nonlinear time-invariant system or a complex stochastic process, and the application of control-theoretical tools to devise intervention policies to stop or mitigate the spreading process. Particular attention should be devoted to the analysis and control of stochastic epidemic models, which are attracting attention in computational epidemiology, physics, and network science, but are still mostly overlooked in the systems and control community.

Besides these research directions, we believe that our scientific community can provide important advances to other open problems of mathematical epidemiology. A key challenge that the scientific community is currently facing with the COVID-19 outbreak is how to integrate the available data within the mathematical models, that is, how to identify the model parameters from clinical and epidemiological data, which are often noisy and incomplete~\cite{Wood2020}. Different optimization methods have been proposed to deal with this problem. See, e.g.,~\cite{Calafiore2020}. From a modeling point of view, an important challenge for the ongoing research, which we briefly mentioned in this survey, is the inclusion of behavioral factors in the model formulation. Game theory has emerged as a promising modeling framework to model human decision-making processes and has been utilized, for instance, to capture the adoption of self-protective behaviors~\cite{Poletti2010,Reluga2010} and to characterize the decision whether to vaccinate or not~\cite{Bauch2004,Trajanovski2018,Hota2019}.

\section*{Acknowledgements}
This work was partially supported by the European Research Council (ERC-CoG-771687) and the Netherlands Organization for Scientific Research (NWO-vidi-14134).

\end{document}